\documentclass{article}
\usepackage[utf8]{inputenc}
\usepackage{amsmath}
\usepackage{amsthm} 
\usepackage{amsfonts} 
\usepackage{amssymb}
\usepackage{verbatim}
\usepackage{mathtools}
    \usepackage[T1]{fontenc}
    \usepackage[utf8]{inputenc}
    \usepackage[british]{babel}
 
 \usepackage[capitalize]{cleveref}
    \usepackage{xcolor}
	\usepackage{soul}

\usepackage[square,numbers,sort&compress]{natbib} %=== sorting references

\usepackage{color}
\definecolor{VeryLightBlue}{rgb}{0.9,0.9,1}
\definecolor{LightBlue}{rgb}{0.8,0.8,1}
\definecolor{MidBlue}{rgb}{0.5,0.5,1}
\definecolor{DarkBlue}{rgb}{0,0,0.6}
\definecolor{Blue}{rgb}{0,0,1}
\definecolor{Gold}{rgb}{1,0.843,0}
\definecolor{LightGreen}{rgb}{0.88,1,0.88}
\definecolor{MidGreen}{rgb}{0.6,1,0.6}
\definecolor{DarkGreen}{rgb}{0,0.6,0}
\definecolor{VeryLightYellow}{rgb}{1,1,0.9}
\definecolor{LightYellow}{rgb}{1,1,0.6}
\definecolor{MidYellow}{rgb}{1,1,0.5}
\definecolor{DarkYellow}{rgb}{1,1,0.01}
\definecolor{DarkPurple}{rgb}{.6,0,1}
\definecolor{Red}{rgb}{1,0,0}
\definecolor{VeryLightRed}{rgb}{1,0.9,0.9}
\definecolor{LightRed}{rgb}{1,0.8,0.8}
\definecolor{MidRed}{rgb}{1,0.55,0.55}

\newcommand{\darkGreen}[1]{{\color{DarkGreen}{#1}}}
\newcommand{\liu}[1]{\darkGreen{\bf [Liu: #1]}}

\newtheorem{theorem}{Theorem}
\newtheorem{defn}{Definition}
    \newtheorem{lem}[theorem]{Lemma}
\newtheorem{conj}{Conjecture}
    \newtheorem{ex}{Example}
    \newtheorem{cor}[theorem]{Corollary}
    \newtheorem{prop}[theorem]{Proposition}

\DeclareMathOperator{\gen}{gen}
\DeclareMathOperator{\Sgen}{S}
\DeclareMathOperator{\Span}{span}
\DeclareMathOperator{\rnk}{\mathnormal{rn_k}}

\DeclareMathOperator{\diam}{diam}
\DeclareMathOperator{\dis}{d}

\newcommand{\Z}{\mathbb{Z}}

\def\leq{\leqslant}
\def\geq{\geqslant}

\def\diam{{\rm diam}}

\def\LB{{\rm LB}}

\DeclarePairedDelimiterX{\norm}[1]{\lVert}{\rVert}{#1}
\DeclarePairedDelimiter{\ceil}{\lceil}{\rceil}

\DeclarePairedDelimiter{\floor}{\lfloor}{\rfloor}

\title{Radio-$k$-Labeling of Cycles for Large $k$\thanks{Research is partially  supported by  the National Science Foundation grant DMS-1600778 and NASA NASA MIRO grant NX15AQ06A.}
}
\author{Colin Bloomfield\thanks{California State University Los Angeles, currently Vanderbilt University; email: \newline colinbloomfield1@gmail.com.}, 
Daphne Der-Fen Liu\thanks{corresponding author.} \  and Jeannette Ramirez\thanks{California State University Los Angeles, Los Angeles, USA; emails: dliu@calstatela.edu  and jeannette.r.831@gmail.com.}
}

\begin{document}
\maketitle

\begin{abstract}
Let $G$ be a simple connected graph. For any two vertices $u$ and $v$, let $d(u,v)$ denote the distance between $u$ and $v$ in $G$. A \textit{radio-$k$-labeling} of $G$ for a fixed positive integer $k$ is a function $f$ which assigns to each vertex a non-negative integer label such that for every two vertices $u$ and $ v $ in $G$, $|f(u)-f(v)| \geq k - d(u,v) +1$. The \textit{span} of $f$ is the difference between the largest and smallest labels of $f(V)$. The  \textit{radio-$k$-number} of a graph $G$, denoted by $rn_k(G)$, is the  smallest span among all radio-$k$-labelings admitted by $G$. A  cycle $C_n$ has diameter $d=\lfloor n/2 \rfloor$. In this paper, we 
combine a lower bound approach with  
cyclic group structure to  
determine  the value of $\rnk(C_n)$ for $k  \geq n-3$. For $d \leq k < n-3$, we  obtain  the values of $\rnk(C_n)$ when $n$ and $k$ have the same parity, and prove partial results when $n$ and $k$ have different parities. 
Our results extend the known values of $rn_d (C_n)$ and $rn_{d+1} (C_n)$ shown by Liu and Zhu \cite{Li05}, and by Karst, Langowitz, Oehrlein and Troxell \cite{Ka17}, respectively.
\end{abstract}
%%%%%%%%%%%


%%%%%%%%%%%
%
\section{Introduction}
%
%%%%%%%%%%%%%%%%%%%%%%

Radio-$ k $-labeling is motivated by the channel assignment problem which was first introduced by Hale \cite{Ha80}. In the channel assignment problem, we must both avoid signal interference between radio stations and minimize the range of frequencies used. Radio-$k$-labeling models this problem with a graph where each station or transmitter is a vertex and neighboring stations are connected by an edge. The distance between two stations is then the usual distance between vertices on a graph, and a station's assigned frequency is represented by a non-negative integer label of the corresponding vertex. If two stations are close enough 
for their signals 
to interfere, they must be assigned frequencies whose difference is determined by their proximity: The closer the stations are geographically, the larger their separation in assigned frequencies must be. 
The goal is to assign frequencies to all stations in a way that minimizes the range of the frequencies used while still avoiding interference. 

Let $G$ be a graph. 
The \textit{distance} between two vertices $u$ and $v$, denoted by $ \dis(u,v)$, is the length of a  shortest path between $u$ and $v$. 
The \textit{diameter} of $ G $ is the largest distance between any two vertices in $V(G)$, which we denote by $\diam(G)$ or sometimes $d$ when the graph is clear in the context. 

	For a positive integer $k$, a function $ f \colon V(G) \to \{0, 1, 2, 3, \cdots\}$ is a \textit{radio}-$k$\textit{-labeling} of a graph $ G$, if for all $u,v \in V(G)$, $u \neq v$, 
	$$
	|f(u)-f(v)| \geq k-\dis(u,v) +1.
	$$
	Sometimes we say $f$ is a \textit{labeling} instead of \textit{radio-$k$-labeling} or omit mentioning the labeled graph $ G $ if it is clear from the context. 
	The \textit{span} of a radio-$k$-labeling $f$ of $G$ is defined as 
	$ \Span(f) \coloneqq \max \{|f(u)-f(v)|: u,v \in V(G) \}.$  
We follow the convention of labeling at least one vertex $ 0 $, so the span of any radio-$k$-labeling is the largest label assigned to a vertex $v \in V(G)$.

The \textit{radio-$k$-number} of $G$ is the minimum span of a radio-$k$-labeling admitted by $G$. 
	We call a radio-$k$-labeling $f$ of a graph $G$ \textit{optimal} if the span is equal to the radio-$k$-number. 
A special case for the  radio-$k$-labeling of a graph $G$ is when $k=\diam(G)=d$, in which a radio-$d$-labeling is also known as {\it radio labeling} and the {\it radio number} of $G$ is denoted by $rn(G)$, where $rn(G)=rn_{d}(G)$.

The notion of radio labeling was  introduced by Chartrand, Erwin and Zhang \cite{Chartrand2}. The radio number for special families of graphs has been studied widely in the literature.  Liu and Zhu  \cite{Li05} determined the radio number of cycles and paths. A general lower bound of the radio number of trees was given by Liu  \cite{Liu}, and then applied by several authors to special families of trees (cf. \cite{Li, Tuza, Bantva_banana_tree, Liu-saha}). Khennoufa and Togni \cite{Khennoufa} studied the radio number for hypercubes by using generalized binary Gray codes. 
Zhou \cite{Zhou1, Zhou2, Zhou3} and Martinez et al. \cite{JP}   investigated the radio number for  Cayley graphs and generalized prism graphs, respectively. 
Niedzialomski \cite{graceful} studied radio graceful graphs (where $G$ admits a surjective  radio labeling from $V(G)$ to $\{1,2,\cdots, |V(G)|\}$) and showed that the Cartesian  product of $t$ copies of a complete graph is radio graceful for certain $t$, providing infinitely many examples of radio graceful graphs of arbitrary diameter.  
For two positive integers $m, n \geqslant 3$, the {\it toroidal grid} $G_{m,n}$ is the Cartesian product of cycles $C_m$ and $C_n$, $C_{m}\Box C_{n}$. Morris et al. \cite{MCC} determined $rn(G_{n,n})$, and Saha and Panigrahi \cite{LP2} determined $rn(G_{m,n})$ when $mn\equiv0\pmod2$. Bantva and Liu \cite{block} considered the radio number for block graphs and line graphs of trees.  

When $k=1$, radio-$k$-labeling is equivalent to vertex coloring. Radio-$k$-labelings of graphs for other values of $k$ have have also been studied  extensively.   There are many known results when $k=2$ (cf. \cite{gri92, cal06}), this case is called $L(2,1)$-labeling, which can be generalized to $L(h,k)$-labeling.  Recently, Chavez et al. studied optimal  radio-$k$-labeling for trees \cite{Chavez-Liu-Shurman}. 

	For cycles $C_n$ with $n \geq 3$, 	the value of $rn_d(C_n)$ was completely determined in \cite{Li05}. 
Also, $rn_{d-1}(C_n)$ when $n \not \equiv 0 \pmod{4} $ has been determined in \cite{Li12}, where bounds are given for the case when $n \equiv 0 \pmod{4} $. Karst, Langowitz, Oehrlein, and Troxell in \cite{Ka17} studied radio-$k$-labeling when $k \geq \diam(C_n)$ and determined the value of $rn_k(C_n)$ when $k=\diam(C_n)+1$.

In this paper we investigate the values of $rn_k(C_n)$ for all  $k \geq   \diam(C_n)$. In particular, we determine the values of $rn_k(C_n)$ for $k \geq n-3$, and for $k < n-3$ with $n \equiv k$ (mod 2).  For $k < n-3$ and $k \not\equiv n$ (mod 2), we establish lower and upper bounds and show these bounds are sharp for some values of $k$ and $n$.  A major tool we use is a combination of a lower bound approach (see \cref{lem: lower bound}) with some basic properties of  
finite cyclic groups. 

%%%%%%%%%%%%%%%%%%5
%
\section{Preliminary and Notations}
%
%%%%%%%%%%%%%%

In what follows, we let $ \mathbb{Z}$ denote the set of integers, $ \mathbb{N}$ the set of natural numbers and $ \mathbb{N}_0 \coloneqq \mathbb{N} \cup \{ 0 \}$. For $n,m \in \mathbb{Z} $, such that $ n \leq m$, we define $ [n,m] \coloneqq \{n, n+1, \ldots, m\}$.  For $n \in \mathbb{N}$, we let $ \mathbb{Z}_n $ denote the cyclic  group of order $n$. For $m \in \mathbb{Z}_n$, let $\langle m \rangle$ be the subgroup of $ \mathbb{Z}_n$ generated by $m$. That is, $\langle m \rangle = \{tm: t \in \Z\}$ (mod $n$). 
Let $ \mathbb{Z}_n / \langle m \rangle $ denote the cyclic quotient group of $ \mathbb{Z}_n$  modulo $ \langle m \rangle$, which is defined on the cosets $ i / \langle m \rangle \coloneqq \{ j \in \mathbb{Z}_n : i-j \in \langle m \rangle \} $ for $ i \in \mathbb{Z}_n  $. 

	We regard $C_n$ as a Cayley graph with vertex set $\Z_n$ and edges generated by $\pm 1$. That is, $V(C_n) = \{0,1,2,\cdots, n-1\}$ where $u$ and $v$ are adjacent if $u \equiv v \pm 1$ (mod $n$). In the proofs we will often make use of an additive structure on $V(C_n)=\Z_n$ which is addition modular $n$. 
	The modest amount of group theory mentioned in this paper can be found in any introductory textbook on abstract algebra. 

By definition, when $k \geq \diam(G)$,  every radio-$k$-labeling of $G$ is injective. So we can order the  vertices of $G$ such that the labeling is strictly increasing. 
Thus, given a radio-$k$-labeling $ f $ of $C_n$ with $k \geq \lfloor n/2 \rfloor$, we denote the vertices of $V(C_n)$, $x_i$ for $i \in [0,n-1]$ such that if $i < j $, then $ f(x_i) < f(x_j) $ for $ i,j \in [0, n-1]. $ 
For $ i \in [0,n-2] $, we define the \textit{label gap} and the \textit{distance gap} respectively by \begin{equation*}
f_i \coloneqq f(x_{i+1}) - f(x_i) 
\ \ \mbox{and} \  \ 
 d_i \coloneqq d(x_{i+1},x_{i}).
\end{equation*}
These notations are used throughout the paper. By definition, $f_i \geq k+1- d_i$. 

\begin{defn}\rm
	When defining a labeling, we first specify an ordering of the vertices $x_0, x_1, \ldots, x_{n-1}$ of $C_n$ by defining a {\it jumping} (or {\it jump})  {\it sequence}  $j_i$ with $1 \leq j_i \leq \frac{n}{2}$, $i \in [0, n-1]$.  This jumping sequence gives the ``clockwise distance" of $x_{i+1}$ from $x_i$, 	%Recall $V(C_n) = \Z_{n}$. 
	namely, if $x_{i}=t \in \Z_n$, then $x_{i+1} = t + j_i$ (mod $n$).
\end{defn}

\begin{defn} \rm
\label{def: per. jump sequence} 
	Let $\mathbf{j} =  (j_{0}, \ldots, j_{t-1}) $, 
	$j_i \in [1,\lfloor n/2 \rfloor]$, $i \in [0,t-1]$.  
	For $n \in \mathbb{N}$, we define the \textit{jump sequence generated by} $ \mathbf{j}$, $\gen_{n}(\mathbf{j}) = (a_{i})_{i=0}^{n-1}$, with $a_{0} = 0$ and $a_{m} = \sum_{i = 0}^{m-1} j_{i}$, where the indices of the jumps are taken modulo $t$ and the sums modulo $n$. We call $ \mathbf{j} $ a {\it periodic jump sequence} with $t$ jumps.
	We call the range of $\gen_{n}(\mathbf{j})$ the set $\textit{generated}$ by $\mathbf{j}$, and denote it as $S_{n}(\mathbf{j})$. 
\end{defn}
		Here is an example of Definition \ref{def: per. jump sequence}. 	
\begin{ex} \rm
\label{ex: per. jumping sequence}
	Let $n=12$ and $\mathbf{j}=(6,4)$, so $t=2$. Then $\gen_{n}(\mathbf{j})= \gen_{12}((6,4)) = (a_0,a_1, \dots$,  
	$a_{11})=(0,6,10,4,8,2,6,0,4,10,2,8).$ And $S_{n}(\mathbf{j})= S_{12}((6,4))=\{0,2,4,6,8,10\}.$ 
\end{ex}
Notice that in \cref{ex: per. jumping sequence}, the jump sequence does not determine a permutation and so does not define a labeling for $C_n$. In practice, we will define a jump sequence and then prove that it determines a permutation of the vertices of $C_n$.

\begin{lem}\label{lem:Partition Size}
	Fix $n \geq 3 $.  Let  $ \mathbf{j} = (j_{0}, j_{1})$,    
		$ h = n-(j_0 + j_1) $, and $ H = \langle h \rangle $. Then
	\begin{equation*}
	|\Sgen_{n}(\mathbf{j})| \leq
	\begin{cases}
	|H| &  \text{if  $j_0 \in H$}, \\
	2|H| & \text{if  $j_0 \not \in H$}.
	\end{cases}
	\end{equation*}
\end{lem}

\begin{proof}
		Note that in $\Z_n$, $H = \langle h \rangle = \langle -h \rangle =\langle j_0+j_1 \rangle$. Hence $j_0+j_1 \in H$. If $j_0 \in H$, then $j_1 \in H$. Thus, $\Sgen_{n}(\mathbf{j}) \subseteq H $, so $ |\Sgen_{n}(\mathbf{j})| \leq |H| $.
	
	If $j_0 \not \in H$, we look at the first $ 2|H| $ vertices obtained from the jump sequence. The set of even and odd   indexed terms, respectively, are  
	\begin{equation*} 
	E = \{ 0,j_0+j_1, 2(j_0+j_1), \ldots,  (|H|-1)(j_0+j_1)\}, \ \ 
	O = j_0 + E. 
	\end{equation*}
	Since $E \subseteq H$ and $j_0 \not \in H$, we have that $E$ and $O$ belong to different cosets of $\mathbb{Z}_n / H$. Hence, $ |O| \leq |H| $ and $ |E| \leq |H|$. Thus, $ |\Sgen_{n}(\mathbf{j})| = |O \cup E| \leq 2|H| $. 
\end{proof}

We will frequently use the following lemma in our proofs.  
\begin{lem}
	\label{lem: gcd(n,h*) when n is even.} 
	Let $n$ be even and $h^* < n$.  Let $\langle h^* \rangle$ be the cyclic subgroup of $\mathbb{Z}_n$ generated by $h^*$. Then   	
	\begin{equation*}
	\gcd(n, h^*) =
	\begin{cases}
	\gcd(\frac{n}{2}, h^*) &  \text{if  $\frac{n}{2} \in \langle h^* \rangle$}, \\
	2 \gcd(\frac{n}{2}, h^*)  & \text{otherwise}.
	\end{cases}
	\end{equation*}
	\end{lem}
\begin{proof}
Let $ n = 2^{l}m $, $ h^{*} = 2^{a}b $, where $ l,m,b \in \mathbb{N} $, $a \in \mathbb{N}_0 $, and both $m$ and $b$ are odd. If $\frac{n}{2} \in \langle h^{*} \rangle $, then $ \gcd(n,h^*) \mid \frac{n}{2} $. It follows that $ a < l $ and so 
$ \gcd(n,h^{*}) = \gcd(\frac{n}{2},h^*)$.
    
    If $ \frac{n}{2} \notin \langle h^{*} \rangle $, then $ \gcd(n,h^{*})  \nmid \frac{n}{2} $ and so $ l \leq a $. It follows that $ \gcd(n,h^{*}) = 2\gcd(\frac{n}{2},h^*) $.
\end{proof}

For any $n \geq 2$ and $k \geq \lfloor n/2 \rfloor$,  the following definition and notation will be used throughout the paper: 
\begin{equation}
\label{phi}
    \varPhi(n,k) \coloneqq \ceil*{\frac{3k+3-n}{2}}. 
\end{equation}
\noindent 
The next result of Karst, Langowitz, Oehrlein, and Troxell \cite{Ka17} plays a critical role in our proofs. For completeness we include a proof.  

\begin{prop} 
\label{lem: lower bound} {\rm \cite{Ka17}} 
For $n > 2$ and $k \geq d$ where $d = \lfloor \frac{n}{2} \rfloor $ and $\varPhi(n,k) \coloneqq \ceil*{\frac{3k+3-n}{2}}$, we have the following lower bounds:
	\begin{equation*}
	\rnk(C_n) \geq
	\begin{cases}
	 \varPhi (n,k) (\frac{n-2}{2}) +k-\frac{n}{2}+1 &\mbox{\rm{if $n$ is even;}} \\
	 \varPhi (n,k)  	(\frac{n-1}{2}) &\mbox{\rm{if $n$ is odd.}} \label{eqn: phi function}
	\end{cases}
	\end{equation*}
\end{prop}

\begin{proof}
	 Assume $f$ is a radio-$k$-labeling of $C_n$. For $i \in [0,n-3]$ the following hold:
	\begin{enumerate}
		\item  $f(x_{i+1})- f(x_i) \geq k+1 - d_i$,
		\item  $f(x_{i+2})-f(x_{i+1}) \geq k + 1 - d_{i+1}$,
		\item  $f(x_{i+2})-f(x_{i}) \geq k+1 - d(x_i,x_{i+2}).$
	\end{enumerate}
	
\noindent
Adding these three inequalities we get 
$$
2(f(x_{i+2})- f(x_i)) \geq 3k+3 - (d(x_i,x_{i+2}) + d_{i+1} + d_i).
$$
	Since $ n \geq d(x_i,x_{i+2}) + d_{i+1} + d_i $, we have
$$
	f_{i+1} + f_i = f(x_{i+2})-f(x_{i}) \geq \ceil*{\frac{3k+3-n}{2}} = \varPhi(n,k).
$$
Hence, we  obtain the lower bounds for the radio-$k$-number of $C_n$:
	\begin{equation*}
	\rnk(C_n) \geq
	\begin{cases}
	\varPhi (n,k) (\frac{n-2}{2}) +k- d+1 & \text{if $n$ is even;} \\
	\varPhi (n,k) (\frac{n-1}{2}) & \text{if $n$ is odd.} \label{eqn: phi function}
	\end{cases}
	\end{equation*}
The result follows. 
\end{proof}

Denote the lower bounds of $rn_k(C_n)$ from 
\cref{lem: lower bound} by ${\rm LB}(n,k)$:  
	\begin{equation}
	\label{lower bound}
	rn_k(C_n) \geq \LB(n,k) := 
	\begin{cases}
	\varPhi (n,k) (\frac{n-2}{2}) +k- \frac{n}{2}+1 &\text{if $n$ is even;} \\
	\varPhi (n,k) (\frac{n-1}{2}) &\text{if $n$ is odd.}
	\end{cases}
	\end{equation}

The next result from \cite{Ka17} greatly reduces the number of inequalities to verify when determining if a map $ f $ is a radio-$k$-labeling.
\begin{prop} {\rm \cite{Ka17}} \label{prop: k < n-3} 
	If $k \geq d$, then $f \colon V(C_n) \to \mathbb{N}_0$ is a radio-$k$-labeling of $C_n$ if and only if the following hold:
	\begin{enumerate}
		\item $f_i \geq k +1 - d_i, \quad \text{for }i \in [0,n-2]$
		\item  $f_i + f_{i+1} \geq k+1 - d(x_i,x_{i+2}), \quad \text{for }i \in [0,n-3]$.
	\end{enumerate}
\end{prop}
\noindent
We do not include a proof of this result, however, the proof of \cref{lem: k >= n-3} in the next section follows similarly.

%%%%%%%%%%%%%%%%%%%%%%%
%%%%%%%%%%%%%%%%%%%%%%%
%
\section{Exact Values for $rn_k (C_n)$}
%
%%%%%%%%%%%%%%%%%%%%%%%
%%%%%%%%%%%%%%%%%%%%%%%

In this section, we completely determine the radio-$k$-number for $C_n$ when  $k \geq n-3$ (\cref{prop: k>=n-3}) and when $k < n-3$ wile $n$ and $k$  have the same parity (\cref{thm: big proof}).

For $k \geq n-3$, we make use of the following lemma, which is a stronger version of Proposition \ref{prop: k < n-3} for the special case when $ k \geq n-3$. 

\begin{lem} 
\label{lem: k >= n-3}
	For $k \geq n-3$, a function $f \colon V(C_n) \to \mathbb{N}_0 $ is a radio-$k$-labeling of $C_n$ if and only if $f_i \geq k - d_i + 1$ for all $i \in [0,n-2]$.
\end{lem}
\begin{proof}
	The forward direction follows from the definition of radio-$k$-labeling. Let $ k \geq n-3 $ and $f$ satisfy $ f_i \geq k - d_i +1$, for all $i \in [0,n-2]$. By Proposition \ref{prop: k < n-3} we only need to show that $f_i + f_{i+1} \geq k+1 - d(x_i,x_{i+2})$ is satisfied.  Let $ i \in [0,n-3]$. Then 
$	f(x_{i+2}) \geq f(x_{i+1}) +k +1 - d_{i+1}$ and 
$	f(x_{i+1}) \geq f(x_i) + k + 1 - d_i.$  
	Combining these two inequalities,
	\begin{equation} \label{eqn: span of 2}
	f(x_{i+2}) \geq f(x_i) + 2k + 2 - (d_{i+1} + d_i).  %\eqno(1)
	\end{equation}
	Since 	$( d_{i+1}+d_{i}) \leq n-1 $ ($\because$ $d_i, d_{i+1} \leq n/2$ and $d_{i+1}+d_{i} \neq n$)  and $k\geq n-3$, we have 
$$
	f(x_{i+2})  \geq f(x_i) + 2k + 2 - (n-1)   \geq f(x_i) + k \geq f(x_i) + k + 1 - d(x_{i},x_{i+2}).    
$$
	Hence the result follows. \end{proof}

\begin{theorem} \label{prop: k>=n-3}
	Let $k \geq n-3 $. Then
	\begin{equation*}
	rn_{k}(C_{n}) = \begin{cases}
	(\frac{n-2}{2})(2k+3-n)+ k - \frac{n}{2} + 1  &  \text{if $n$ is even}; \\
	(\frac{n-1}{2})(2k+3-n)  &  \text{if $n$ is odd}.
	\end{cases}
	\end{equation*}
\end{theorem}

\begin{proof}
Suppose $f$ is a radio-$k$-labeling for $C_n$. From 
Equation (\ref{eqn: span of 2}):
	\begin{equation*}
	f({x_{i+2}}) - f(x_{i}) \geq 2k+2-(d_{i}+d_{i+1}) \geq 2k+3-n.
	\end{equation*}
	Let $\Phi^{\prime}(n,k) = 2k+3-n $. 
	If $n$ is odd, then
	\begin{align*}
	\Span(f) \geq \sum_{i = 0}^{(n-3)/2}[f(x_{2i+2})-f(x_{2i})] \geq \left(\frac{n-1}{2}\right)\Phi^{\prime}(n,k). 	\end{align*}
	If $n$ is even, the last vertex cannot be paired with others and so,	
	$$
	\Span(f) \geq \sum_{i = 0}^{(n-4)/2}[f(x_{2i+2})-f(x_{2i})] + f_{n-2}  \geq \left(\frac{n-2}{2}\right)\Phi^{\prime}(n,k)+k+1-n/2.
	$$
Therefore, the lower bounds are established. 

		Next we find a labeling which meets the lower bounds. 
Suppose $n = 2q+1 $ for $ q \in \mathbb{N}$. Let the distance gap  be $d_i = q$ 
for all $i \in [0,n-2]$. Since $ d = q $ and $ \gcd(2q+1,q) = 1 $, if we let $j_i = d_i $, then $j_i$ determines a permutation of $\mathbb{Z}_n $. 
Define the label gap  by $ f_i = k + 1 - d$. Then $ f_i \geq k+1-d_i$ for all $i \in [0,n-2]$.  Therefore, $f$ is a valid radio-$k$-labeling of $C_n$ and $ \Span(f)$ meets the desired lower bound.

	Now consider the case when $n$ is even. That is, $n=2q$ for $ q \in \mathbb{N}$, and $ d = q $. For $i \in [0,n-2]$ we define the jumping and labeling sequences  respectively by:
	\begin{equation*}
	j_i = \begin{cases}
	d &\text{$i$ is even}; \\ 
	d-1 &\text{$i$ is odd}.
	\end{cases}
	\hspace{0.3in}
	f_i = \begin{cases}
	k+1-d &\text{$i$ is even}; \\
	k+2-d &\text{$i$ is odd}.
	\end{cases}\hspace{0.2in}
	\end{equation*}
	It is easy to see that the jump sequence gives a permutation of $V(C_n)$, and it holds that $ f_{i} = k+1 -d_{i} $ for all $i$. %If $ i $ is odd, $ f_{i} = k+1 -(d-1) = k+1-d_{i} $. 
	By \cref{lem: k >= n-3}, $f$ is a radio-$k$-labeling of $C_n$ 
	with the desired span.
\end{proof}

%%%%%%%%%%%%%%%%%%%%%%%%%%%%%%%%%%%
%
Now we turn to the case when $d \leq k<n-3$ and $n  \equiv k$ (mod $2$).  
Recall  $\varPhi(n,k)$ and $\LB(n,k)$ from \cref{phi} and \cref{lower bound}, respectively.  

\begin{lem}
\label{two values of h}
	Let $f$ be a radio-$k$-labeling of $C_n$, where $n$ and $k$ have the same parity.  For any $i \in [0,n-3]$, if $f_i+f_{i+1}=\varPhi(n,k)$, then $d(x_i,x_{i+2}) \in \{ \frac{n-k}{2},\frac{n-k}{2}-1 \}$.
\end{lem}

%%%%%%%%%%%%%

\begin{proof}
	Suppose $f_i+f_{i+1}=\varPhi(n,k)$.  Because $n$ and $k$ have the same parity, 
	by the definition of radio-$k$-labeling, we have 
$$
(3k-n+4)/2 = 	\varPhi(n,k) = f_i+f_{i+1}  \geq 2(k+1) - (d_i + d_{i+1}).
$$
Hence, $d_i + d_{i+1} \geq \frac{n+k}{2}$. Combining this with the fact that $ d(x_i,x_{i+2}) \leq n- (d_i + d_{i+1}) $ we get $d(x_{i}, x_{i+2}) \leq (n-k)/2$. Moreover, because  $f(x_{i+2})-f(x_i)=f_i+f_{i+1}
=\varPhi(n,k) \geq k+1- d(x_i,x_{i+2})$, we obtain   $ d(x_i,x_{i+2}) \geq 
(n-k-2)/2$. Therefore, the result follows. 
\end{proof}

\cref{two values of h} shows that when $n$ and $k$ have the same parity, in order to keep the $\varPhi$-function value for any three consecutive vertices, $x_i, x_{i+1}, x_{i+2}$, the distance between $x_i$ and $x_{i+2}$ must be either $\frac{n-k}{2}$ or $\frac{n-k}{2}-1$.

\begin{theorem}
	\label{thm: big proof} If $ n $ and $ k $ have the same parity and $ \lfloor \frac{n}{2} \rfloor \leq k \leq n-3 $, then $rn_k(C_n)=\LB(n,k)$.
\end{theorem}

\begin{proof}
	Throughout the proof (in fact, throughout the paper), in every defined labeling, $d< d_i + d_{i+1} \leq n$ holds for all $ i $ considered, which implies that $n- (d_i + d_{i+1})= d(x_i, x_{i+2}).$ 
	Also, in every defined jump sequence,
	no jump $j_i $ will exceed the diameter of the graph, so that $j_i = d_i.$
	
	\noindent \textbf{Case 1.} Both $ n $ and $ k $ are even. 
	Denote $ n = 4q + t$ and $ k = 2q + 2m + t $, where $ t \in \{0,2\}$ and $ 0 \leq m < q-1 $. 
	Let $ h = q-m-1 $ and $ p = \gcd(d,h) $. Define the jump sequence 
	(or  distance gap) by:

	\begin{equation*}
	j_{i} = d_i = \begin{cases}
	d & i  \text{ is even,} \\
	d-h &  i \text{ is odd and } i \neq \frac{nx}{p} -1, \ x \in \mathbb{N}, \\
	d-h-1 & i \text{ is odd and } i= \frac{nx}{p} -1 \ \mbox{for some} \ x \in \mathbb{N}.
	\end{cases}
	\end{equation*}
	
	\noindent \textbf{Claim:} The jump  sequence gives a permutation. 

	\noindent
	{\it Proof of Claim)} Recall from Section 1  that $ V(C_{n} = \mathbb{Z}_{n} $. If $ p = 1 $, then our jumping strategy goes as follows:  Start at the vertex $x_0$, then  jump $d$ to the next vertex $x_1$. The vertices $x_0$ and $x_1$ form a coset of $ \mathbb{Z}_{n}/\langle d \rangle.$ (Note, each coset in $\mathbb{Z}_{n}/\langle d \rangle$ consists of two elements, $\{i, i+d\}$ (mod $n$).) After $x_1$ we then jump $d-h$ to $x_2$, which is in another coset of $\mathbb{Z}_{n}/\langle d \rangle$. After $x_2$ we then jump $d$ to the vertex $x_3$. The vertices $x_2$ and $x_3$ form another  coset of $\mathbb{Z}_{n}/\langle d \rangle$. 
	
		We repeat this process until each vertex is landed exactly once. This can be done since 
	$\gcd(d,h)= \gcd(d, d-h)=1$. So $|S_{n}(\mathbf{j})|=n$ with $\mathbf{j}=(d,d-h)$. Therefore, $\mathbf{j}=(d,d-h)$ is a permutation on $\Z_n.$
	
	If $ p \neq 1 $, by Lemma \ref{lem: gcd(n,h*) when n is even.},  
$\gcd(n,h)= \gcd(d,h) = p$ if $d \in \langle h \rangle$; otherwise 
$\gcd(n,h)=2 \gcd(d,h) = 2p$.  
	If $d \in \langle h \rangle$,  then the order of each coset in $\Z_n/ \langle h \rangle$ is $n/p$.   Also,  $S_{n}(\mathbf{j})$ with  $\mathbf{j}=(d,d-h)$ covers one coset of $\Z_n/\langle h \rangle$.  So, $|S_{n}(\mathbf{j})|=n/p$.
	Similarly, if $d \not \in \langle h \rangle$, by \cref{lem:Partition Size},    $S_{n}(\mathbf{j})$ with 
	$\mathbf{j}=(d,d-h)$ covers two cosets in $\Z_n/\langle h \rangle$.  Thus, $|S_{n}(\mathbf{j})|=n/p$.
	
	After making $\frac{n}{p}-1$ jumps, we jump $d-h-1$ to the next vertex, which is $-1 + \langle h \rangle$ in $ \mathbb{Z}_{n}/ \langle h \rangle $. Then we continue jumping $\mathbf{j}=(d,d-h)$, which will cover one or two new cosets, depending on whether $d \in \langle h \rangle$, as shown above. We continue this process until all vertices are covered exactly once. Thus the jump sequence determines a permutation. 
	\hfill$\blacksquare$

	%%%%%%%%%%%%%%%%%%%%%%%%%%%
	%%%%%%%%%%%%%%%%%%%%%%%%%%%%%%%		
	
	Next we define the labeling sequence as follows: 
	\begin{equation*}
	f_i= 
	\begin{cases}
	2m+1 + (t/2) & \quad \text{if $i$ is even}, \\
		q+m+1 + (t/2) & \quad \text{if $i$ is odd}.
	\end{cases}
	\end{equation*}
To complete the proof for Case 1, it suffices to show: 

\noindent	
\textbf{Claim:} $f_i$ is a valid radio-$k$-labeling with the desired span. 	
	
\noindent
{\it Proof of Claim)}	By Proposition \ref{prop: k < n-3}, we need to show that $f_i \geq k+1-d_i$ and $f_i + f_{i+1} \geq k+1-d(x_i, x_{i+2}).$ We call these the first and the second inequality, respectively.
	We begin by verifying the first inequality. 
	Suppose $i$ is even. Then $j_i = d_i = d = 2q + (t/2)$  and 
	\begin{equation*}
	k+1-d_i= 2q+2m+t+1-2q-(t/2)=2m+1+(t/2)  = f_i. %\liu{replace by =}
 	\end{equation*}
	Suppose $i$ is odd. Then $j_i = d_i \in \{q+m+1+(t/2), q+m+(t/2) \}$ and 
	$$ 
	k+1-d_i \leq 2q+2m+t+1-(q+m+t/2)= q+m+1 + (t/2)= f_i.
	$$
	Now we verify the second inequality. By the facts that $d(x_i,  x_{i+2})=n-(d_i+d_{i+1})$  and 
	$d_i+d_{i+1} \in \{3q+m+t, 3q+m+t+1\}$, we obtain   
	\begin{align*}
	f_i+ f_{i+1} 
	 = q+3m+t+2 & = 3q+m+t+1- 2q+2m+1 \\
	& \geq d_i+d_{i+1} - 2q+2m+1 \\
		& = k+1 -  d(x_i,x_{i+2}). 	
	\end{align*}
	Therefore, $f$ is a  radio-$k$-labeling with the   desired span.  \hfill$\blacksquare$

\begin{comment}
      	If $ n = 4q + 2 $ and $k=2q+2m+2$, then the labeling sequence is defined by: 
	\begin{equation*}
	f_{i} = \begin{cases}
	2m+2 & \text{if $ i $ is even}, \\
	q+m+2 & \text{if $ i $ is odd}.
	\end{cases}
	\end{equation*}
	To show that $f$ is a radio-$k$-labeling follows similarly to the $n = 4q$ case. 	For an example of Case 1, 
	%where $n$ and $k$ are both even 
	see Example \ref{ex: case1} following the proof.
	\liu{end}
	\end{comment}
	
	\begin{comment}
	\noindent \textbf{Claim:} $f_i$ is a valid radio-$k$-labeling.  We begin by verifying the first inequality.
	
	\noindent	Suppose that $i$ is even, then
	\begin{equation*}
	k+1-d_i= 2q+2m+3-2q-1=2m+2 \leq f_i.
	\end{equation*}
	Suppose $i$ is odd, then $j_i \in \{q+m+2, q+m+1 \}$ and
	$$ k+1-d_i \leq 2q+2m+3-(q+m+1)=q+m+2 = f_i.$$
	
	\noindent	Next we verify the second inequality. We use the fact that $d(x_i, x_{i+2})=n-(d_i+d_{i+1})$ and $d_i+d_{i+1} \in \{3q+m+2,  3q+m+3 \}.$ So we have
	\begin{align*}
	f_i+ f_{i+1} & = q+3m+4 \notag \\
	& = 3q+m+3- (2q+1-2m)+2 \\
	& \geq d_i+d_{i+1} - (2q+1-2m)+2 \\
	& = 2q+1-4q-2+2m +d_i+d_{i+1} +2 \\
	& = 2q+1+2m-[4q+2 -d_i+d_{i+1}] +2 \\
	& = k- d(x_i, x_{i+2})+1.
	\end{align*}
	Therefore, we have a valid radio-$k$-labeling with span: $ \Span(f)= \frac{n-2}{2}(f_i + f_{i+1})+f_{n-2}.$ We have $f_i + f_{i+1}=q+3m+4= \varPhi(n,k)$ and $f_{n-2}=2m+2=k+1-d$. So $ \Span(f)= \frac{n-2}{2}\varPhi(n,k)+k+1-d$, which is exactly the lower bound in Proposition \ref{lem: lower bound}.
\end{comment}
\medskip

		\noindent \textbf{Case 2.} Both $ n $ and $ k $ are odd. Let $ n = 4q+3 $ and $ k = 2q + 2m + 1 $ for $ 0 \leq m \leq q-1 $, or $ n = 4q+1 $ and $ k = 2q+2m+1 $ for $ 0 \leq m \leq q-2 $. Then $\varPhi(n,k)= (3k-n+4)/2$.  Let 
	$$ z = \ceil*{\frac{d+q+m+1}{2}}  \hspace{0.2in}   \mbox{and} \hspace{0.2in}  
	s = \frac{n}{\gcd(n,z)}.
	\hspace{0.5in}
	$$
	
%%%%%%%%%%%%%%%
\begin{comment}
	
	Then
	\begin{equation*}
	\varPhi(n,k)=
	\begin{cases}
	q+3m+2 &  \text{if  $n = 4q+3$}, \\
	q+3m+3 & \text{if  $n = 4q+1$}.
	\end{cases}
	\end{equation*}
	%%%%%%%%%%%%%%%%%%%%%%%%%%%%%%%%	
	\begin{equation*}
	z =
	\begin{cases}
	\ceil*{\frac{3q+m+2}{2}} &  \text{if  $n = 4q+3$}, \\
	\ceil*{\frac{3q+m+1}{2}} & \text{if  $n = 4q+1$}.
	\end{cases}
	\end{equation*}

\end{comment}
%%%%%%%%%%%%%%%	
%	\liu{Replace Sub-Case}
	
	\noindent 
	\textbf{Case 2.1.} Suppose $\varPhi(n,k)=(3k-n+4)/2$ is even. Define the jump sequence and labeling sequence respectively by:
	\begin{equation*}
	j_i=
	\begin{cases}
	z +1 &  \text{if  $i = sx-1$}, \ x \in \mathbb{N} \\
	z & \text{otherwise}; 
	\end{cases} 
	\hspace{0.2in} 
\mbox{and} 
	\hspace{0.2in} 
	f_i=
	\frac{\varPhi(n,k)}{2} \quad \text{for all $i$}.
	\end{equation*} 
We claim that our jump sequence gives a permutation.	Let $ V(C_{n}) = \mathbb{Z}_{n} $. If gcd$(n,z)=1$, then by the second formula of $j_i$ in the above, $ \langle z \rangle = \mathbb{Z}_{n}$ gives a permutation. If gcd$(n,z) > 1$,  since $ s = |\langle z \rangle| $, the jump sequence labels each coset of $ \mathbb{Z}_{n}/ \langle z \rangle $, and then  transitions into the next unlabeled coset by jumping $z+1$ (the first formula in $j_i$ above).  Continue this process until all vertices are labeled. Hence, the jump sequence is a permutation.

It remains to show that $f_i$ is a valid radio-$k$-labeling. Consider the case that $n=4q+3$ and $k=2q+2m+1$. 
	Then   
	$\varPhi(n,k)= 
	q+3m+2$. By the assumption that $\varPhi(n,k)$ is even,  
	$q$ and $m$ must have the same parity, and 
	\begin{equation*}
	j_i=
	\begin{cases}
	\frac{3q+m+4}{2} &  \text{if  $i = sx-1$}, \ x \in \mathbb{N} \\
	\frac{3q+m+2}{2} & \text{otherwise};
	\end{cases} 
	\hspace{0.2in} 
	\mbox{and} 
	\hspace{0.2in} 
	f_i=
	\frac{q+3m+2}{2} \quad \text{for all $i$}.
	\end{equation*} 
	To verify the first inequality, 
as $d_i = j_i \geq (3q+m+2)/2$, we have 
	\begin{align*}
	k+1-d_i \leq 2q+2m+2- \frac{3q+m+2}{2} = \frac{q+3m+2}{2} =f_i.
	\end{align*}
	To verify the second inequality, using the facts that $d_i+d_{i+1}= n - d(x_i, x_{i+2})$ and  $d_i+d_{i+1} \in \{3q+m+2, 3q+m+3\}$, we obtain  
	\begin{align*}
	f_i+ f_{i+1}  = q+3m+2 
	& = 3q+m+3- (2q-2m+1) \\
	& \geq d_i+d_{i+1} - (2q-2m+1) \\
	& = k- d(x_i, x_{i+2})+1.
	\end{align*}
	Therefore, $f$ is a  radio-$k$-labeling with the desired span. 
	The case when $n=4q+1$ follows similarly.

%%%%%%%%%%%%%%%%%%%%	
	\begin{comment}
	Next, consider $ n = 4q+1 $ and $k=2q+2m+1$.  Then $\varPhi(n,k)= \ceil*{ \frac{3k-n+3}{2}} = q+3m+3$ and $z= \ceil*{ \frac{3q+m+1}{2}}$. Since $q+3m+3$ is even, $q$ and $m$ have different parities,  it follows that
	%%%%%%%%%%%%%%%%%%%%%%%%%%%%%%	
	\begin{equation*}
	j_i=
	\begin{cases}
	\frac{3q+m+3}{2} &  \text{if  $i = sx-1$}, \ x \in \mathbb{N} \\
	\frac{3q+m+1}{2} & \text{otherwise},
	\end{cases}
	\end{equation*}
	%%%%%%%%%%%%%%%%%%%%%%%%	
	\begin{equation*}
	f_i=
	\frac{q+3m+3}{2} \quad \text{for all $i$}.
	\end{equation*}
	
	We begin by showing that the first inequality holds. Since $d_i \geq \frac{3q+m+1}{2}$, then
	\begin{align*}
	k+1-d_i \leq 2q+2m+2- \frac{3q+m+1}{2} = \frac{q+3m+3}{2} = f_i.
	\end{align*}
	Next we verify the second inequality; notice that $f_i+f_{i+1}=q+3m+3$ and $d_i+d_{i+1} \in \{3q+m+1,3q+m+2\}.$ Hence,
	%%%%%%%%%%%%%%%%%%%	
	\begin{align*}
	f_i+ f_{i+1} & = q+3m+3 \notag \\
	& = 3q+m+2- (2q-2m)+1 \\
	& \geq d_i+d_{i+1} - (2q-2m)+1 \\
	& = 2q-4q+1+2m +d_i+d_{i+1}  \\
	& = 2q -4q-1+2+2m+d_i+d_{i+1}  \\
	& = 2q+2m+2-[n-(d_i+d_{i+1} )] \\
	& \geq k- d(x_i, x_{i+2})+1.
	\end{align*}
	Therefore, we have a valid radio-$k$-labeling with $ \Span(f)= \frac{n-1}{2}(f_i + f_{i+1}).$ We have $f_i + f_{i+1}=q+3m+3= \varPhi(n,k)$. So $ \Span(f)= \frac{n-1}{2}\varPhi(n,k)$, which is exactly the lower bound in Proposition \ref{lem: lower bound}. 
	\end{comment}
%%%%%%%%%%%%%%%%%%%%%

	\noindent \textbf{Case 2.2.} Suppose $\varPhi(n,k) = (3k-n+4)/2 $ is odd. If $ n = 4q+3 $, $\varPhi(n,k) = q + 3m +2 $, so $ q $ and $m$ must have different parities. It follows that $d+q+m+1 = 3q+m+2$ is odd. Otherwise, $n = 4q +1 $ and $ \Phi(n,k) = q + 3m+3 $. Then $q$ and $m$ have the same parity and again $ d+q+m+1 $ is odd. Hence $z=(d+q+m+2)/2$ always holds.  	Let $h=n-2z$. 

%%%%%%%%%%%%%%%%%%%%%%%%%%
	\begin{comment}
	Then  %Then $z= \lceil \frac{d+q+m+1}{2} \rceil = \frac{d+q+m+2}{2}.$  
	\begin{equation*}
	z =
	\begin{cases}
	\frac{3q+m+3}{2} &  \text{if  $n = 4q+3$}, \\
	\frac{3q+m+2}{2} & \text{if  $n = 4q+1$}.
	\end{cases}
	\end{equation*}
	
	We have
	\begin{equation*}
	h =
	\begin{cases}
	q-m &  \text{if  $n = 4q+3$}, \\
	q-m-1 & \text{if  $n = 4q+1$}.
	\end{cases}
	\end{equation*}
	\end{comment}
	
%%%%%%%%%%%%%%%%%%%

	\noindent \textbf{Case 2.2.1.} Suppose $\gcd(n,h)=1.$ Define the jump sequence and the labeling sequence respectively by:
	$$ 
	j_i = z \quad \text{for all $i$}
	\hspace{0.4in} 
\mbox{and} 
\hspace{0.4in}
	f_i=
	\begin{cases}
	\floor{ \frac{\varPhi(n,k)}{2}} &  \text{if  $i$ is even},  \\
	\ceil{\frac{\varPhi(n,k)}{2} }  & \text{if  $i$ is odd}.
	\end{cases}
	$$ 
To show that the jump sequence is a permutation, we need to verify that $\langle z \rangle = \Z_n$. Since $\gcd(n,h)=1$, we have that $\langle h \rangle = \mathbb{Z}_n = \langle -h \rangle = \langle 2z \rangle $. This implies that $\gcd(n,2z)=1$. Hence $\gcd(n,z)=1$ and so $\langle z \rangle = \Z_n$, implying the jump sequence is a permutation. 

	Now we show $f_i$ is a valid radio-$k$-labeling. We prove the case when $n=4q+1$. (The case when $n=4q+3$ follows similarly.) Then  $\varPhi(n,k)=q+3m+3$.  
%%%%%%%%%%%5
\begin{comment}	
	So:
	\begin{equation*}
	f_i=
	\begin{cases}
	\lfloor\frac{q+3m+3}{2} \rfloor &  \text{if  $i$ is even},  \\
	\lceil\frac{q+3m+3}{2} \rceil  & \text{if  $i$ is odd}.
	\end{cases}
	\end{equation*}	
\end{comment}
%%%%%%%%%%%%	
	For the first inequality, one can see that  
	$ k+1- d_i = %2q+2m+2- \frac{3q+m+2}{2} = \frac{q+3m+2}{2}= 
	\floor*{\frac{\varPhi(n,k)}{2} } \leq f_i$. 
	For the second inequality, since $d_i+d_{i+1}= n - d(x_i, x_{i+2})$ and $d_i+d_{i+1}=3q+m+2$, one can verify that $f_i+ f_{i+1} \geq 
	k- d(x_i, x_{i+2})+1$.
		
		Therefore, we have a valid radio-$k$-labeling with the desired span,  which is exactly the lower bound in Proposition \ref{lem: lower bound}. 
	
%%%%%%%%%%%%%%%%%%%%%%%%%%%%%%%%%%%%%%%%%%%%%%%%%
%
\noindent \textbf{Case 2.2.2.}  \label{case: 2.2.2.}
	Suppose $c=\gcd(n,h)>1$ for some odd $c$. 
	Then $\langle c \rangle = \langle h \rangle$, and there are $c$ cosets in $\mathbb{Z}_n / \langle c \rangle$. %\liu{end} 
	Recall that 
	\begin{equation*}
	h =
	\begin{cases}
	q-m &  \text{if  $n = 4q+3$}, \\
	q-m-1 & \text{if  $n = 4q+1$}.
	\end{cases}
	\end{equation*}
	
	Define the jumping  sequence and labeling sequence respectively by:
	\begin{equation*}
	j_i=
	\begin{cases}
d & \quad \text{if $i$ is even, $i \in [0,n-2-\frac{n}{c}]$}, \\
		d -h  & \quad \text{if $i$ is odd, $ i = \frac{2n}{c}x-1, x \in \mathbb{N}$, $i \in [0,n-1-\frac{n}{c}]$},
		\\
d -h +1 & \quad \text{if $ i $ is odd, $ i \not= \frac{2n}{c}x-1, \ x \in \mathbb{N}$, and $i \in [0,n-1-\frac{n}{c}]$, } \\
	(n-h)/2 & \quad \text{if $i \in [n-\frac{n}{c},n-2]$};
	\end{cases}
	\end{equation*}

	\begin{equation*}
	f_i=
	\begin{cases}
	k+1 - d & \quad \text{if $i$ is even and $i \in [0,n-2-\frac{n}{c}]$}, \\
	
	\varPhi(n,k)- (k+1 - d) & \quad \text{if $ i $ is odd and $i \in [0,n-1-\frac{n}{c}]$}, \\
	
	\lfloor \frac{\varPhi(n,k)}{2} \rfloor & \quad \text{if $i$ is even and $i \in (n-1-\frac{n}{c},n-2]$}, \\

	\lceil \frac{\varPhi(n,k)}{2} \rceil  & \quad \text{if $i$ is odd and $i \in (n-\frac{n}{c},n-2]$}.
	\end{cases}
	\end{equation*}
	
Consider the case that $n=4q+1$ and $k=2q+2m+1$. (The case for $n=4q+3$ and $k=2q+2m+1$ follows  similarly.) Note that $d=2q$,  $\varPhi(n,k)=q+3m+3$, and $h=q-m-1$. 

%%%%%%%%%%%%%
	\begin{comment}

	Using the above, our jumping sequence and labeling sequence are:
	\begin{equation*}
	j_i=
	\begin{cases}
	2q & \quad \text{if $i$ is even and $i \in [0,n-2-\frac{n}{c}]$}, \\
	
	q+m+1  & \quad \text{if $ i = \frac{2n}{c}x-1, \quad x \in \mathbb{N}  $ and $i \in [0,n-2-\frac{n}{c}]$}, \\
	
	q+m+2 & \quad \text{if $ i $ is odd and $ i \not= \frac{2n}{c}x-1, \quad x \in \mathbb{N}$ and and $i \in [0,n-2-\frac{n}{c}]$, } \\
	\frac{3q+m+2}{2} & \quad \text{if $i \in (n-2-\frac{n}{c},n-2]$};
	\end{cases}
	\end{equation*} 
	
	\begin{equation*}
	f_i=
	\begin{cases}
	2m+2 & \quad \text{if $i$ is even and $i \in [0,n-2-\frac{n}{c}]$}, \\
	
	q+m+1 & \quad \text{if $ i $ is odd and $i \in [0,n-2-\frac{n}{c}]$}, \\
	
	\frac{q+3m+2}{2} & \quad \text{if $i$ is even and $i \in (n-2-\frac{n}{c},n-2]$}, \\

	\frac{q+3m+4}{2}  & \quad \text{if $i$ is odd and $i \in (n-2-\frac{n}{c},n-2]$}.
	\end{cases}
	\end{equation*}
	\end{comment}
%%%%%%%%%%%%%%%

\noindent \textbf{Claim:} Our jumping sequence is a permutation.  

	\noindent
	{\it Proof of Claim)} 
	We show that we label two cosets of $\mathbb{Z}_n / \langle c \rangle$ at a time by  alternately jumping $2q$ and $q+m+2$ (see the first and third formulas in $j_i$ above). That is, $\mathbf{j}=(2q, q+m+2)$. 
	
Note that $h = n - (3q+m+2)$. This implies $\langle h \rangle = \langle 3q+m+2 \rangle = \langle c \rangle$. Hence, it is enough to show that $2q \not\in \langle c \rangle$.	
Because $\gcd(4q+1, 2q)=1$, we have $ \langle 2q \rangle = \Z_n$. 
If $2q \in \langle c  \rangle$, then $\langle c \rangle = \Z_n$, contradicting to gcd$(n,h)=c > 1$. Thus, $2q \notin \langle c \rangle$, and our jump sequence begins by labeling two cosets in $\mathbb{Z}_n / \langle c \rangle$.

	Next we show that the first two cosets  labeled  are $\langle c \rangle$ and $\langle c \rangle + \lfloor \frac{c}{2} \rfloor$, by proving that $2q \in \lfloor \frac{c}{2} \rfloor + \langle c \rangle$. Let $n=4q+1 = tc$ for some odd $t$. Then 
$$
2q = \frac{tc-1}{2} = \frac{c(t-1)}{2} + \frac{c-1}{2}.  
$$
Thus, $2q \in \lfloor \frac{c}{2} \rfloor + \langle c \rangle$, and so $2q + \langle c \rangle=  \lfloor \frac{c}{2} \rfloor + \langle c \rangle$. Hence, the first two labeled cosets are $ \langle c \rangle$ and $\langle c \rangle + \lfloor \frac{c}{2} \rfloor $. 

Afterwards, we jump $q+m+1$ (the second formula in $j_i$ given above) and land in the coset $\langle c \rangle -1$ in $\Z_n / \langle c \rangle$. Continue the same process repeatedly as above, we label the first $c -1 $ cosets in pairs in the following order: 
	%%%%%%%%%%%%%%%%%%%
$$ 
\left( \langle c \rangle, \langle c \rangle + \floor*{\frac{c}{2}} \right) \rightarrow \left(\langle c \rangle -1, \langle c \rangle + \floor*{\frac{c}{2}} - 1 \right) \rightarrow \cdots \rightarrow \left(\langle c \rangle +1 - \floor*{\frac{c}{2}}, \langle c \rangle + 1 \right).
$$
	%%%%%%%%%%%%%%%%%%%%
After the last pair of cosets above, the only coset that was left unlabelled is $ \langle c \rangle + \lceil \frac{c}{2} \rceil$. Again, we jump 		$q+m+1$ and land in that coset. 

It remains to show that jumping $\frac{n-h}{2}$ (the last formula in the $j_i$ given above) is a  permutation of $ \langle c \rangle + \lceil \frac{c}{2} \rceil$.  
This is true because $n$ is odd, so $\gcd(n,\frac{n-h}{2})=\gcd(n,n-h) = \gcd(n,h)=c$.  
\hfill$\blacksquare$

Finally we show that $f_i$ is a radio-$k$-labeling. 
	The first inequality can be obtained directly from our definitions that for all $i$, 
	$k+1-d_i \leq f_i.$ 
	Next we verify the second inequality by listing all possible cases in Table \ref{table:2.2}.   Note, for the last case (last row), we need to show $\frac{3q+5m+4}{2} \geq \frac{q+7m+6}{2}$.  Recall $0 \leq m \leq q-2$. Then $q \geq m+1$, which implies $\frac{3q+5m+4}{2} \geq \frac{q+7m+6}{2}$. For an example of this case, see Example \ref{ex: Case 2.2.2.}.
	
	Therefore, we have a valid radio-$k$-labeling with the desired span. The proof is complete.
\end{proof}

\begin{center}
	\begin{table}[t]
		\begin{tabular}{||c || c || c || c || c||}
			\hline
			$i \in$ & $i+1 \in$ & $d_i + d_{i+1} \leq$ & $f_i + f_{i+1}$ & $k-d(x_i +x_{i+2})+1 \leq$\\ [0.5ex] 
			\hline\hline
			$[0, n-1-  \frac{n}{c}]$ & $[0, n-1- \frac{n}{c}]$ & $3q+m+2$ & $q+3m+3$ & $q+3m+3$\\ [1.5ex]
			& 
			&  &  & \\ [1.5ex]
			\hline 
			$[n-\frac{n}{c}, n-2]$ & $[n-\frac{n}{c}, n-2]$ &  $3q+m+2$ & $q+3m+3$ & $q+3m+3$ \\ [1.5ex]
	
			 &  &  &  & \\ [1.5ex]
			\hline
			$ \{ n-1-\frac{n}{c} \} $ & $ \{ n - \frac{n}{c} \}  $ & $\frac{5q+3m+4}{2}$ & $\frac{3q+5m+4}{2}$ & $\frac{q+7m+6}{2}$ \\ [2ex]
			
			$i$ is odd & 
			$i$ is even &  &  & \\ [1.5ex]
			\hline
			
		\end{tabular}
		\caption{All possibilities for checking the second inequality for Case 2.2.2.}
%		\colin{We fixed the issues with the chart.}
		\label{table:2.2}
	\end{table}
\end{center}

	\begin{ex} \rm
	\label{ex: case1} 
		In Figure \ref{fig:n=16 one coset} we let $q=4$ and $m=1$, so $n=16$ and $k=10$. We get $\varPhi(16,10)=9$, $h=2$, and $p= \gcd(d,h)=2$. In this example we label one coset of $\mathbb{Z}_n / \langle h \rangle$ at a time. With $\mathbf{j}=(8,6)$, $S_{16}((8,6))$ covers one coset (circular black vertices) of $\Z_{16} / \langle 2 \rangle$. We then jump $5$ into a new coset (gray diamond vertices) of $\Z_{16} / \langle 2 \rangle$ to complete the labeling. 
		
		\begin{figure}[h]
			\centering
			\includegraphics[width=0.63\linewidth]{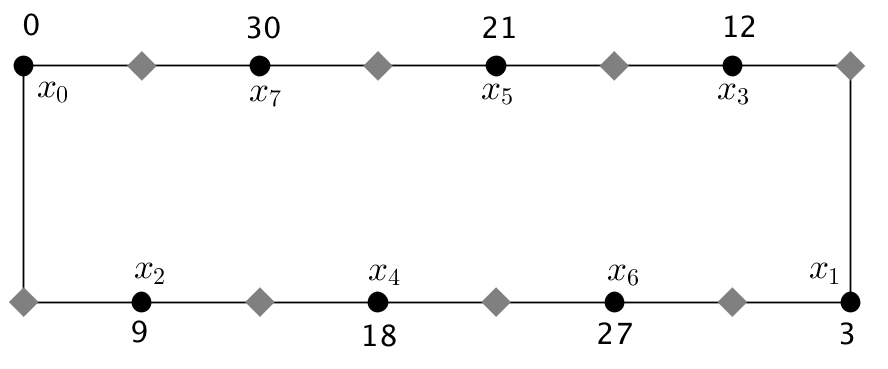}
			
			\vspace{10px}
			
			\includegraphics[width=0.63\linewidth]{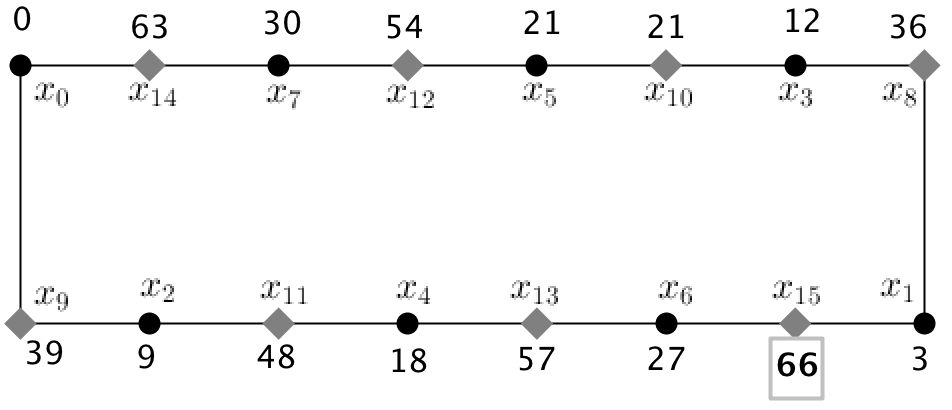}
			
			\vspace{10px}

			\caption{Steps of finding an optimal radio-$10$-labeling for $C_{16}$ in Example \ref{ex: case1}.}
			
			\label{fig:n=16 one coset}
		\end{figure}
	\end{ex}
	
	\begin{ex} \rm
\label{ex: Case 2.2.2.}
	In Figure \ref{fig:n=27 three cosets} we let  $n=27$ and $k=19$. Then $d=13$ and $h=3$. Notice that $\gcd(n,h)= \gcd(27,3)=3=h$ and  $\varPhi(27,19)=17$. 
	We label two cosets of $\mathbb{Z}_n / \langle h \rangle$ at a time. With $\mathbf{j}=(13,11)$, $S_{27}((13,11))$ covers two cosets (solid black and white circles) of $\Z_{27} / \langle 3 \rangle$. We then jump $10$ into a new coset (gray diamonds) of $\Z_{27} / \langle 3 \rangle$. We then jump $12$ consecutively to finish labeling the last coset.
	
	\begin{figure}[h]
		\centering
		\includegraphics[width=1\linewidth]{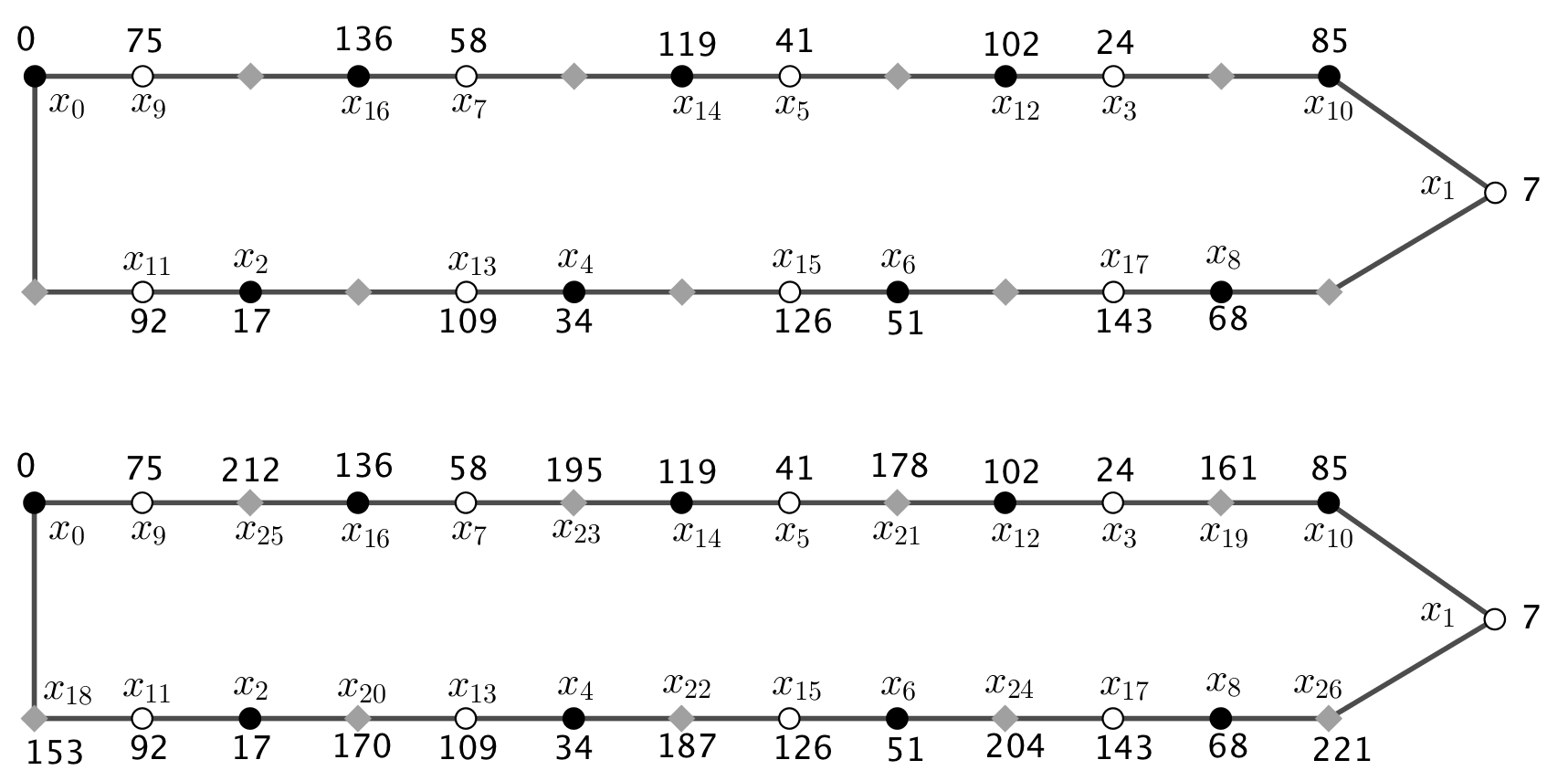}
		
		\caption{Example \ref{ex: Case 2.2.2.}: Steps of finding an optimal radio-$19$-labeling for $C_{27}.$}
		\label{fig:n=27 three cosets}
	\end{figure}

\end{ex}

%%%%%%%%%%%%%%%
%%
\section{Results for $k<n-3$ and $n \not \equiv k$ (mod $2$)  }
%%
%%%%%%%%%%%%%%%

The case when $n$ and $k$ have different parities is more complicated, and so more interesting!  Following the work in the previous section when $n$ and $k$ have the same parity, we will use the lower bound proved in 
\cref{lem: lower bound}. On the other hand, there are properties we use in the previous section that are not valid anymore for the case when $n$ and $k$ have different parities. The following lemma reveals the fact that in searching for an  optimal radio-$k$-labeling that achieves the lower bound for $C_n$ we would lose some flexibility 
 we had in \cref{two values of h}. 
 Consequently,  we show that the lower bound in \cref{lem: lower bound} is not always achieved.    
 
   \begin{lem}\label{lem: To obtain phi sum of two jumps must be -h.}
   \label{fix}
   	Let $f$ be a radio-$k$-labeling for $C_n$, where $n$ and $k$ have different parities.  For any $i \in [0,n-3]$, if $f_i+f_{i+1}=\varPhi(n,k)=  (3k-n+3)/2$, then $d(x_i,x_{i+2})=k+1-\varPhi(n,k)= (n-k-1)/2$.
   \end{lem}

\begin{proof}
This lemma essentially follows from Lemma 2.1 in \cite{Ka17}. One may also prove this lemma with the same approach taken in the proof of Lemma \ref{two values of h}.
\end{proof}

%%%%%%%%%%%%%%%%%%%
\begin{comment}
  \begin{proof}
	Suppose $f_i+f_{i+1}=\varPhi(n,k)$.  Because $n$ and $k$ have different parity, 
% 	$\varPhi(n,k)= \frac{3k-n+4}{2}.$
	by the definition of radio-$k$-labeling, we have 
$$
(3k-n+3)/2 = 	\varPhi(n,k) = f_i+f_{i+1}  \geq 2(k+1) - (d_i + d_{i+1}).
$$
Hence, $d_i + d_{i+1} \geq \frac{n+k+1}{2}$. Combining this with the fact that $ d(x_i,x_{i+2}) \leq n- (d_i + d_{i+1}) $ we get $d(x_{i}, x_{i+2}) \leq (n-k-1)/2$. Moreover, because  $f(x_{i+2})-f(x_i)=f_i+f_{i+1}
=\varPhi(n,k) \geq k+1- d(x_i,x_{i+2})$, we obtain   $ d(x_i,x_{i+2}) \geq 
(n-k-1)/2$. Therefore, the result follows. 
\end{proof}

   \begin{proof}
   	Suppose $n \not\equiv k$   (mod 2), and   $f_i+f_{i+1}=\varPhi(n,k)=(3k-n+3)/2$. 
  Because $f_i \geq k+1-d_i$ and 
   	$f_{i+1} \geq k+1 - d_{i+1}$,  
we have 
$   	\varPhi(n,k)=f_i+f_{i+1} \geq 2(k+1) - (d_i + d_{i+1}).$ 
Thus, 
$d_i + d_{i+1} \geq (k+n+1)/2.$ 

On the other hand, by definition, 
%   	\begin{align*}
 $  	f(x_{i+2})-f(x_i)=f_i+f_{i+1}\geq k+1- d(x_i,x_{i+2})$. Hence,  %$\frac{3k-n+3}{2} \geq  k+1- d(x_i,x_{i+2}) \\
  % 	\Rightarrow d(x_i,x_{i+2}) & \geq k+1- \frac{3k-n+3}{2} \\
  % 	\Rightarrow
$d(x_i,x_{i+2})  \geq (n-k-1)/2$.
%   	Using $\frac{n-k-1}{2} \leq d(x_i,x_{i+2})$ with the 
Because $d(x_i,x_{i+2})=n-(d_{i}+d_{i+1})$, we get 
   	$d_{i+1}+d_{i} \leq (k+n+1)/2$.  
Hence, $d_{i+1}+d_{i}=(k+n+1)/2$, so  $d(x_i, x_{i+1}) = (n-k-1)/2$.  \end{proof}
\end{comment}

%%%%%%%%%%%%%%
   
  \cref{fix} shows that when $n$ and $k$ have different parities, in order to keep the $\varPhi$-function value for any consecutive three vertices, $x_i, x_{i+1}, x_{i+2}$, the distance between $x_i$ and $x_{i+2}$ must be fixed as   $\frac{n-k-1}{2}$. Throughout the section we denote this fixed value by $h$. That is,  
   \begin{equation}\label{eqn: definition of h}
   h \coloneqq k+1-\varPhi(n,k) =\frac{n-k-1}{2}.
   \end{equation}
%
%%%%%%%%%%%%%%%
%%%%%%%%%%%%%%%
\begin{comment}
\begin{ex}
	Let $n=60$ and $h^*=9$. We get $\frac{n}{2}=30 \in \langle 9 \rangle$. And $\gcd(n,h^*)=\gcd(60,9)=3$ and $\gcd(d,h^*)=\gcd(30,9)=3.$ So $\gcd(n,h^*)=\gcd(d,h^*).$ If $h^*=8$, then $\frac{n}{2}=30 \not \in \langle 8 \rangle$. And $\gcd(n,h^*)=\gcd(60,8)=4$ and $\gcd(d,h^*)=\gcd(30,8)=2.$ So $\gcd(n,h^*)=2\gcd(d,h^*).$
\end{ex}
\end{comment}
%%%%%%%%%%%%%%
%%%%%%%%%%%%%%
% 
This observation gives us the following lower bound. Recall the definition of $\LB(n,k)$ from \cref{lower bound}.
\begin{prop}
\label{prop: General lower bound for n and k different parity.}
    Suppose $ k < n-3 $ and $ n \not \equiv k \pmod{2}$, and let $ p = \gcd(n,h)$, where $h$ is defined in \cref{eqn: definition of h}. Then
    \begin{equation}
        \LB(n,k) + \left\lceil \frac{p}{2} \right\rceil - 1 \leq rn_k(C_n).  
    \end{equation}
\end{prop}
\begin{proof}
    The proof follows from  \cref{lem:Partition Size} and \cref{lem: To obtain phi sum of two jumps must be -h.} which collectively imply that one can label at most $(2n)/p $ vertices while obtaining the $\Phi$-function value.
\end{proof}

When $n$ is even and $k$ is odd, we establish lower and upper bounds as follows and show the bounds are sharp for some $k$ and $n$. 

\begin{theorem}
     \label{thrm: odd bounds}
	Suppose $n$ is even and $k$ is odd with $ \frac{n}{2} \leq k < n-3$. 
Let	$p=\gcd(n,h)$. 
	Then 
	\begin{equation*}
	\LB(n,k) + \ceil*{ \frac{p}{2}} -1 \ \leq rn_k(C_n) \ \leq \LB(n,k) + p -1 . 
	\end{equation*}
	Moreover, if $d \notin \langle h \rangle$ then the lower bound in the above is sharp.    	
\end{theorem}

\begin{proof}
	The lower bound follows from Proposition \ref{prop: General lower bound for n and k different parity.}. As $d= \frac{n}{2}$, to prove the upper bound and the ``moreover'' part, by Lemma \ref{lem: gcd(n,h*) when n is even.}, it suffices to define a radio-$k$-labeling $f$ which satisfies all the following:   
	\begin{itemize}
		\item If $d \in \langle h \rangle$, then $f$ labels one coset of $\mathbb{Z}_n/ \langle h \rangle$ at a time and increases the span by 1 each time it switches to an unlabeled coset.   
		\item If $d \not\in \langle h \rangle$, then $f$ labels two cosets of $\mathbb{Z}_n/ \langle h \rangle$ at a time and increases the span by 1 each time it switches to a new pair of cosets. 	
	\end{itemize}	
Let $p^* = \gcd(d, h)$. 
As $n$ is even and $k$ is odd,  
we let $n=4q$ and $n=4q+2$; and in either case, let $k=2q+2m+1$. Define the jumping sequence and labeling sequence (for both cases) respectively by:
	\begin{equation*}\label{eq:n=4q k odd gap seq}
	j_i=
	\begin{cases}
	d & \quad \text{if $i$ is even}, \\
	
	d-h-1 & \quad \text{if $ i = \frac{n}{p^*}x-1, \quad 0<x<p^*  $}, \\
	
	d-h & \quad \text{if $ i $ is odd and $ i \not= \frac{n}{p^*}x-1, \quad 0<x<p^* $};
	\end{cases}
	\end{equation*} 	
	\begin{equation*}\label{eq:n=4q k odd gap seq}
	f_i=
	\begin{cases}
	k+1-d & \quad \text{if $i$ is even}, \\
	
	d-h+1 & \quad \text{if $ i = \frac{n}{p^*}x-1, \quad 0<x<p^*  $}, \\
	
	d-h & \quad \text{if $ i $ is odd and $ i \not= \frac{n}{p^*}x-1, \quad 0<x<p^* $}.
	\end{cases}
	\end{equation*}   		
\noindent	
Then $f$ defines a radio-$k$-labeling; the proof follows similarly to the proofs in Theorem \ref{thm: big proof}. Whenever $j_i = d-h-1=q+m $, we have $ f_{i-1} +f_i= \varPhi(n,k)+1$. This happens $p^*-1$ times while labeling. Thus, %\liu{Replace "the labeling $f$ gives us a span of"} 
	span$(f)=\LB(n,k) + p^*-1$. 
	
	If $d \in \langle h \rangle$, then $\gcd(n,h) = \gcd(d,h)= p = p^*$, so the upper bound holds. If $d \not\in \langle h \rangle$, then $\gcd(n,h) = 2 \gcd(d,h)$, or $p = 2p^*$;  hence the lower bound is sharp.     
\end{proof}

\begin{cor} \label{cor: n even k odd}
	Suppose $n$ is even and $k$ is odd with $ \frac{n}{2} \leq k < n-3$. 
	Recall that $h= \frac{n-k-1}{2}$. If $\gcd(\frac{n}{2},h) =1$, then $rn_k(C_n)=\LB(n,k)$. 
\end{cor}
\begin{proof}
	Assume $d=\frac{n}{2} \in \langle h \rangle$.  Then by our assumption, $p= \gcd(n,h)=\gcd(\frac{n}{2},h)=1$. From Theorem \ref{thrm: odd bounds} we get that the upper and lower bounds meet and $rn_k(C_n)=\LB(n,k).$ If  $d \not \in \langle h \rangle$, we have $p= \gcd(n,h)=2\gcd(d,h)=2$. The moreover part in Theorem \ref{thrm: odd bounds} implies that 
	$rn_k(C_n)=\LB(n,k) + \lceil \frac{p}{2} \rceil -1 = \LB(n,k)$.
\end{proof}

%%%%%%%%%%%%
%%%%%%%%%%%%
\begin{comment}
  			\begin{ex} \label{ex: n even k odd}
  				
  				In Figure \ref{fig:n=20 one coset} we let  $n=20$ and $k=13$. We get $h=3$ and $\varPhi(20,13)=11$. Since $\gcd(d,h)= \gcd(10,3)=1$, with $\mathbf{j}=(10,7)$ we have $S_{20}((10,7))= \Z_{20}$. And $f_i=4$ if $i$ is even and $f_i=7$ if $i$ is odd. Then $\LB(n,k)=11 \cdot 9 +13-10+1=103$, that is , the lower bound is achieved.  
  				
  				\begin{figure}[h]
  					\centering
  					\includegraphics[width=0.9\linewidth]{cycle20.png}
  					
  					\vspace{10px}

  					\caption{Example \ref{ex: n even k odd}: An optimal radio $13$-labeling for $C_{20}$.}
  					
  					\label{fig:n=20 one coset}
  				\end{figure}
  			\end{ex}
 \end{comment} 	
 %%%%%%%%%%%%%%%
 %%%%%%%%%%%%%%%

We conjecture that the upper bound in \cref{thrm: odd bounds} is tight  when $d \in \langle h \rangle$ (see Conjecture \ref{conj}).

Now we turn to the case that $n$ is odd and $k$ is even. The remaining two results of this section are devoted to this case. Recall $\varPhi(n,k)$, $\LB(C_n)$, and $h$  defined in \cref{phi}, \cref{lower bound}, and \cref{eqn: definition of h}, respectively. 

\begin{theorem}\label{prop: n odd k even}
	Let $n$ be odd and $k$ even. 
	%Denote  $h=k+1- \varPhi(n,k)$. 
	If $\gcd(n,h)=1$ and $h$ is odd, then $rn_k(C_n)=\LB(n,k).$
	\end{theorem}

	\begin{proof}
		We define a jump sequence and labeling as follows, for all $i$:
		$$	
		j_i= \frac{n-h}{2}  \text{ \quad and \quad} 	f_i=\frac{\varPhi(n,k)}{2}.
		$$
It is easy to see that our jumping sequence gives a permutation, since  $\gcd(n, \frac{n-h}{2})=\gcd(n,h)=1$, thus  $\langle \frac{n-h}{2} \rangle= \mathbb{Z}_n.$ 

Now we show that $f$ is a valid radio-$k$-labeling by verifying 1 and 2 in \cref{prop: k < n-3}. 
The first inequality is easy to see. 
				For the second inequality, since $d(x_i,x_{i+2})= n - (d_i + d_{i+1})=h$, we obtain 
		$  
		k + 1 - d(x_i,x_{i+2}) = k +1 -h = %k+1-k-1 + \varPhi(n,k)=
		\varPhi(n,k)= f_i + f_{i+1}.
		$ 
	Therefore, $f$ is a radio-$k$-labeling with the desired span as in \cref{lem: lower bound}.
\end{proof}

\begin{theorem}\label{thm: n odd k even new jump}\label{thm: n odd k even new jump}
Suppose $n=4q+i$ and $k=2q+2m+i-1$ where $i \in \{1,3\}$ and $ 0 \leq m \leq q-2$. 
	If $h \mid n$, then	
	\begin{equation*}
	{rn}_k(C_n) = \LB(n,k) + \frac{h-1}{2}.
	\end{equation*}
\end{theorem}

\begin{proof}
Recall from \cref{eqn: definition of h}, $h = (n-k-1)/2$. Since $n=4q+i$ and $k=2q+2m+i-1$ we get $h = q - m$. 
	That $\LB(n,k) + \frac{h-1}{2} \leq rn_k(C_n)$ holds follows from 
	\cref{prop: General lower bound for n and k different parity.}. 
	
To show that $rn_k(C_n) \leq \LB(n,k) + \frac{h-1}{2}$, we define a jumping and a labeling sequence  in the following:
	\begin{equation*}
	j_i=
	\begin{cases}
	d & \text{if $i$ is even and $i \in [0,n-2-\frac{n}{h}]$}, \\
	
	d-h & \text{if $i$ is odd, $ i = \frac{2n}{h}x-1, x \in \mathbb{N}$, $i \in [0,n-\frac{n}{h}]$}, \\
	
	d-h+1 & \text{if $ i $ is odd and $ i \not= \frac{2n}{h}x-1, x \in \mathbb{N}$, $i \in [0,n-2-\frac{n}{h}]$, } \\
	\frac{n-h}{2} & \text{if $i \in [n-\frac{n}{h},n-2]$};
	\end{cases}
	\end{equation*}	
	\begin{equation*}\label{eq:n=4q k odd gap seq}
	f_i=
	\begin{cases}
	k+1-d & \text{if $i$ is even and $i \in [0,n-2-\frac{n}{h}]$}, \\
	
	\varPhi(n,k)-(k+1-d)+1 & \text{if $i$ is odd, $ i = \frac{2n}{h}x-1, x \in \mathbb{N}$, $i \in [0,n-\frac{n}{h}]$}, \\
	
	\varPhi(n,k)-(k+1-d) & \text{if $ i $ is odd and $ i \not= \frac{2n}{h}x-1, x \in \mathbb{N}$, $i \in [0,n-2-\frac{n}{h}]$, } \\
	\frac{\varPhi(n,k)}{2}  & \text{if $i \in [n-\frac{n}{h},n-2]$}.
	\end{cases}
	\end{equation*} 
	Finishing the proof by showing our jump sequence is a permutation and our labeling sequence is valid follows similarly to the proof of Theorem \ref{thm: big proof}. \end{proof}

Below is an example of Theorem \ref{thm: n odd k even new jump}.

\begin{ex} \rm
\label{ex: thm 4.9}
	In Figure \ref{fig:n=27 two cosets} we let  $n=27$ and $k=20$. Then $d=13$ and $h=3$. Notice that $\gcd(n,h)= \gcd(27,3)=3=h$ and  $\varPhi(27,20)=18$. 
	We label two cosets of $\mathbb{Z}_n / \langle h \rangle$ at a time. With $\mathbf{j}=(13,11)$, $S_{27}((13,11))$ covers two cosets (solid black circles and gray diamonds) of $\Z_{27} / \langle 3 \rangle$. We then jump $10$ into a new coset (larger gray circles) of $\Z_{27} / \langle 3 \rangle$. Notice that the difference of $f(x_{18})-f(x_{16})=163-144=19= \varPhi(27,20)+1$.
	
	\begin{figure}[h]
		\centering
		\includegraphics[width=1\linewidth]{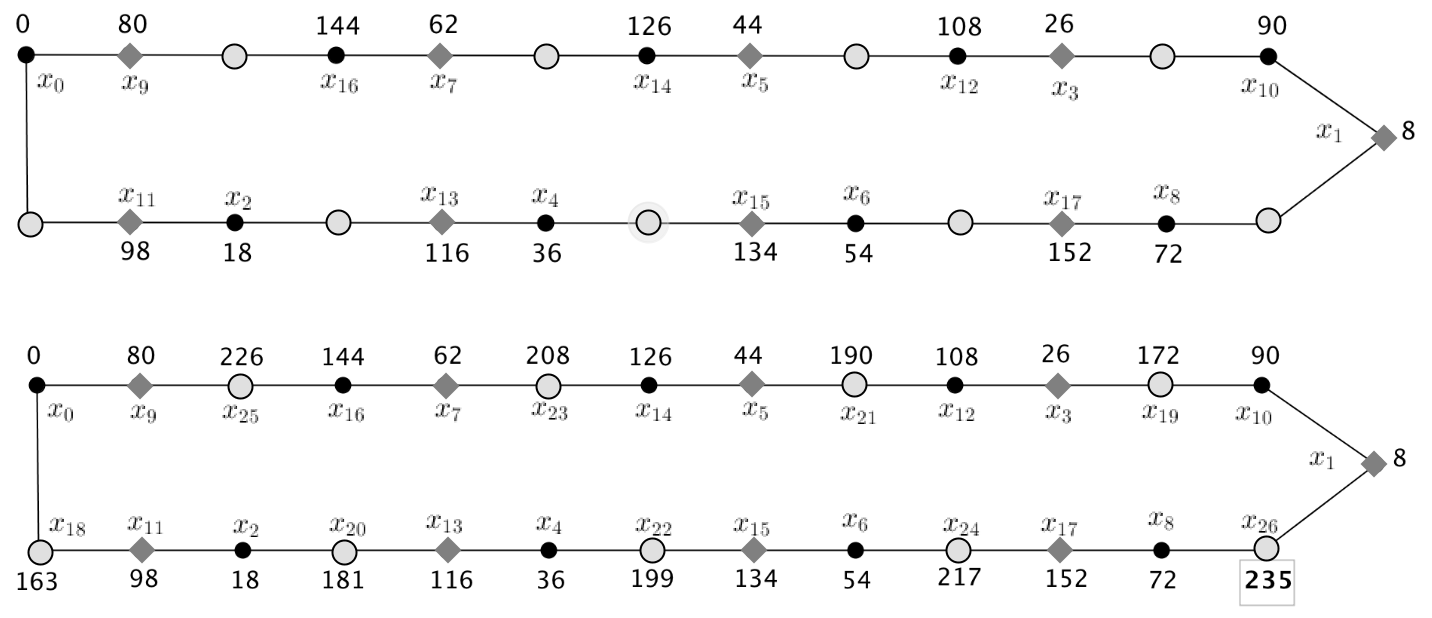}
		
		\caption{Example \ref{ex: thm 4.9}: Steps of finding an optimal radio $20$-labeling for $C_{27}.$}
		\label{fig:n=27 two cosets}
	\end{figure}

\end{ex}

  	\section{Closing Remarks }
  	
   Most results in this paper  are extensions of results by Karst, Langowitz, Oehrlein, and Troxell \cite{Ka17}, and by Liu and Zhu \cite{Li05}. 
  	   	The method of using the  $\varPhi$-function to prove a lower bound for the radio $k$-number for $C_n$ was introduced in \cite{Li05}  when $k =\floor*{\frac{n}{2}}$, the diameter of $C_n$.  
%  	   	\liu{modified this sentence}  
Precisely, it was shown that $rn_{\lfloor n/2 \rfloor}  (C_n)=\LB(n, d)$ for all $n$. 
  	  	This method was applied and    extended in 
  	%by Karst, Langowitz, Oehrlein, and Troxell 
  	\cite{Ka17}, where the authors completely determined the radio $k$-number for $C_n$ when $k = \lfloor \frac{n}{2} \rfloor + 1$, the diameter of $C_n$ plus one. The lower bound is shown again to be close to the exact value. 
  	\begin{theorem} 
  	\label{thm: ka}
{\rm \cite{Ka17}}
  		Let $n=4q+r$ with $q \geq 0$ and $0 \leq r \leq 3.$ For $k = \lfloor \frac{n}{2} \rfloor + 1$,
  		\begin{align*} 
  		&\text{{\rm (i)}} \ r=0: \  rn_{k}(C_n)=
  		\begin{cases}
  		\LB(n,k) & \text{if $q$ is even};\\
  		\LB(n,k)+1 & \text{if $q$ is odd.}
  		\end{cases} \\
  		&\text{{\rm (ii)}} \ r=1,2: \  rn_{k}(C_n)= \LB(n,k).\\
  		& \text{{\rm (iii)}} \ r=3: \  rn_{k}
  		(C_n)=
  		\begin{cases} 
  		\LB(n,k) & \text{if $q \neq 2$ and it not a multiple of 3};\\
  		\LB(n,k)+1 & \text{otherwise.}
  		\end{cases}
  		\end{align*}
  	\end{theorem}
  
%  	This paper extends the results from \cite{Ka17}. 
%The results in this paper agree with the results in \cite{Li05, Ka17}. 
%  	Consider Theorem \ref{thm: liu}: 
%In the case when $n$ and $k=\lfloor \frac{n}{2} \rfloor$ have the same parity, that is $n=4q$ with $k=2q$ and $n=4q+3$ with $k=2q+1$, Theorem \ref{thm: big proof} determines the radio $\lfloor \frac{n}{2} \rfloor$-number which is $LB(n,k).$ 
 
% \liu{modify below} \colin{We think the writing is good, but if you want to make changes feel free.}
Extending this lower bound approach with properties of cyclic groups, we were able to  determine the values of  $rn_k(C_n)$ for most  $k \geq \lfloor n/2 \rfloor$. Many results in  \cite{Li05, Ka17} can be immediately obtained by this paper. For instance,   
  	when $n$ and $d=\lfloor \frac{n}{2} \rfloor$ have different parities we have $n=4q+2$ with $k=2q+1$, or  $n=4q+1$ with $k=2q$. 
  	%Recall from Chapter 3 that $h=\frac{n-k-1}{2}.$ 
  	When $n=4q+2$ and $k=2q+1,$ we have $h=q$. Thus, $\gcd(d,h)= \gcd(2q+1,q)=1.$ Using Corollary \ref{cor: n even k odd}, we have $rn_d(C_n)=\LB(n,d)$. When $n=4q+1$ and $k=2q$ we have $h=q$. Thus $\gcd(n,h)= \gcd(4q+1,q)=1$. If $h$ is odd, by Proposition \ref{prop: n odd k even}, $rn_d(C_n) = \LB(n,d)$. 
%  	Our results agree with Theorem \ref{thm: ka} in the case that $n$ and $k=\lfloor \frac{n}{2} \rfloor +1$ have the same parity and some cases when $n$ and $k$ have different parities. When $n$ and $k$ have the same parity, that is $n=4q+1$ with $k=2q+1$ and $n=4q+2$ with $k=2q+2$, Theorem \ref{thm: big proof} determines the radio $(\lfloor \frac{n}{2} \rfloor +1)$-number which is $LB(n,k).$

For the case that $k=d+1$,   	  	when $n$ and $d+1=\lfloor \frac{n}{2} \rfloor +1$ have different parities we have $n=4q$ with $k=2q+1$, or $n=4q+3$ with $k=2q+2$. %Recall that $h=\frac{n-k-1}{2}.$
  	  	When $n=4q$ and $k=2q+1$, we have $h=q-1$.  Consider the case when $q$ is even. Then,  $\gcd(d,h)=\gcd(2q,q-1)=\gcd(q+1,q-1)=1$.  Using Corollary \ref{cor: n even k odd}, we have 
  	  	$rn_{d+1}(C_n) = \LB(n, d+1)$. 
%  	  	that the radio $(\lfloor \frac{n}{2} \rfloor +1)$-number is equal to $\LB(n,k).$ 
When $n=4q+3$ and $k=2q+2$ we have $h=q$. Thus, $\gcd(n,h)= \gcd(4q+3,q)=\gcd(3,q)=1$ if $q$ is not a multiple of $3$. If $h=q$ is odd, by Proposition \ref{prop: n odd k even} the radio $k$-number is $\LB(n,k).$
  	
Although the values of $rn_k(C_n)$ for most $k$ and $n$ with $n \geq \lfloor n/2 \rfloor$ have been determined, it remains open for other cases, especially when $k$ and $n$ have different parities. %As suggested in our work, utilizing the lower bound method with further group structure is one approach for further investigation.  
So far the evidence we obtained  indicates that the upper bound in \cref{thrm: odd bounds} might be always sharp for $n/2 \in \langle h \rangle $.   
  	\begin{conj}
\label{conj}
  		Suppose $n$ is even and $k$ is odd with $ \frac{n}{2} \leq k < n-3$. Let $h= \frac{n-k-1}{2}$ and $p=\gcd(n,h)$. If $ n/2 \in \langle h \rangle $, then
  		\begin{equation*}
  		rn_k(C_n) = \LB(n,k) + p -1.
  		\end{equation*}
  	\end{conj}
  	
\bigskip

\noindent
{\bf Acknowledgement.} The authors are grateful to the two anonymous  referees for their careful reading of the manuscript and for their constructive suggestions which led to a better presentation of the article.

%With that in mind, we would like to further investigation on \cref{conj} might .    
\begin{comment}
  	
 	\section{Conjecture}
  Recall the general upper and lower bounds from Theorem \ref{thrm: odd bounds} when $n$ is even and $k$ is odd. We conjecture that the upper bound is always tight.
  	\begin{conj}
  		Suppose $n$ is even and $k$ is odd with $ \frac{n}{2} \leq k < n-3$. Let $h= \frac{n-k-1}{2}$ and $p=\gcd(n,h)$. If $ d \in \langle h \rangle $, then
  		\begin{equation*}
  		rn_k(C_n) = \LB(n,k) + p -1.
  		\end{equation*}
  		
  	\end{conj}

  	The conjecture states that when $d \in \langle h \rangle$, the best strategy is to label one coset of $\Z_n /  \langle h \rangle$ at a time. In the following we show that our conjecture matches the result found in \cite{Ka17} when $n$ and $k$ have different parities.

  	When $n=4q$, $k=2q+1$, and $q$ is odd, the authors of \cite{Ka17}  proved that $rn_k(C_n)=\LB(n,k)+1$, as shown in Theorem \ref{thm: ka} . When $n=4q$ and $k=2q+1$,  then $h= \frac{n-k-1}{2}=q-1$. Recall from Observation \ref{obs} that when $d \in \langle h \rangle$, then $\gcd(n,h)=\gcd(d,h)$.
  	We calculate $p=\gcd(n,h)= \gcd(d,h)= \gcd(2q,q-1)=\gcd(q+1,q-1)$, and since $q$ is odd, $p=2$. By our conjecture $rn_k(C_n)=\LB(n,k)+1$ when $k= \lfloor \frac{n}{2} \rfloor +1$. This matches the result of Theorem \ref{thm: ka}. 
  	
\end{comment}

\bibliographystyle{plain}
\bibliography{main.bib}

\end{document}